**Michiel Hazewinkel**　　　　　1　　　　　CWI
Direct line: +31-20-5924204　　　　POBox 94079
Secretary: +31-20-5924233　　　　1090GB  Amsterdam
Fax: +31-20-5924166
E-mail: mich@cwi.nl

# Symmetric functions, noncommutative symmetric functions and quasisymmetric functions II

by

*Michiel Hazewinkel*
*CWI*
*POBox 94079*
*1090GB  Amsterdam*
*The Netherlands*

**Abstract**. Like its precursor this paper is concerned with the Hopf algebra of noncommutative symmetric functions and its graded dual, the Hopf algebra of quasisymmetric functions. It complements and extends the previous paper but is also selfcontained. Here we concentrate on explicit descriptions (constructions) of a basis of the Lie algebra of primitives of *NSymm* and an explicit free polynomial basis of *QSymm*. As before everything is done over the integers. As applications the matter of the existence of suitable analogues of Frobenius and Verschiebung morphisms is discussed.

## 1. Introduction

As said before, [24], the symmetric functions are an exceedingly fascinating object of study; they are best studied from the Hopf algebraic point of view (in my opinion), although they carry quite a good deal more important structures, indeed so much that whole books do not suffice, but see [26, 27, 31, 33, 34].

The first of the two generalizations to be discussed is the Hopf algebra, *NSymm*, of *noncommutative symmetric functions* (over the integers). As an algebra, more precisely a ring, this is simply the free associative ring over the integers, $\mathbf{Z}$, in countably many indeterminates

$$NSymm = \mathbf{Z}\langle Z_1, Z_2, \cdots\rangle \tag{1.1}$$

and the coalgebra structure is given by the comultiplication determined by

$$\mu : Z_n \mapsto \sum_{i+j=n} Z_i \otimes Z_j, \quad \text{where} \quad Z_0 = 1 \tag{1.2}$$

and $i$ and $j$ are in $\mathbf{N} \cup \{0\} = \{0, 1, 2, \cdots\}$. The augmentation is given by

$$\varepsilon(Z_n) = 0, \quad n = 1, 2, 3, \cdots \tag{1.3}$$

(and, of course $\varepsilon(Z_0) = \varepsilon(1) = 1$). The Hopf algebra *NSymm* is a noncommutative covering generalization of the Hopf algebra of symmetric functions,

$$Symm = \mathbf{Z}[z_1, z_2, \cdots] \tag{1.4}$$

where the $z_n$ are seen as either the elementary symmetric functions $e_n$ or the complete symetric functions $h_n$. The interpretation of the $z_n$ as the $h_n$ seems to work out somewhat nicer, for instance in obtaining the standard inner product autoduality of *Symm* in terms of the natural duality between *NSymm* and *QSymm*, the Hopf algebra of quasisymmetric functions, see [24], section 6. *QSymm* will be described and discussed later in this paper.

The projection is given by

$$NSymm \longrightarrow Symm, \quad Z_n \mapsto z_n \tag{1.5}$$

and is a morphism of Hopf algebras.

The systematic investigation of *NSymm* as a noncommutative generalization of *Symm* was started in [14] and continued in a whole slew of subsequent papers, e.g. [7, 8, 9, 20, 21, 22, 23, 25, 28, 29, 30, 32, 46].

It is amazing how much of the theory of *Symm* has natural noncommutative analogues. This includes Newton primitives, Schur functions, representation theoretic interpretations, determinental formulas now involving the quasideterminants of Gel'fand - Retakh, [12, 13]), Capelli and Sylvester identities, and much more. And, not rarely, the noncommutative versions are more elegant than their commutative counterparts.

Note, however, that in most of these papers the noncommutative symmetric functions are studied over a fixed field $K$ of characteristic zero and not over the integers (or a field of positive characteristic). This makes quite a difference, see section 3 below. The papers [19, 20, 21, 22, 23] focuss on the case over the integers, as does the present paper.

It should be stressed that *NSymm* attracts a lot of attention not only as a natural generalization of *Symm*. It turns up spontaneously. For instance in terms of representations of the Hecke algebras at zero, [8, 24, 30, 46] and as the direct sum of the Solomon descent algebras of the symmetric groups, [1, 10, 14, 35, 43, 44] and [39], Ch. 9. Moreover there are e.g. applications to noncommutative continued fractions, Padé approximants, and a variety of interrelations with quantum groups and quantum enveloping algebras, [2, 14, 29, 37]. Further, the duals, the quasisymmetric functions, first turned up (under that name) in the theory of plane partitions and counting permutations with given descent sets, [15, 16, 45]. Actually, *QSymm*, precisely as the graded dual of *NSymm*, goes back at least to 1972 in the theory of noncommutative formal groups, [5]. See [20] for an outline of the role played by *QSymm* in that context. An application of *NSymm* to chromatic polynomials is in [11].

Given a Hopf algebra $H$, with multiplication $m$ and comultiplication $\mu$, a primitive in $H$ is an element $P$ of $H$ such that

$$\mu(P) = 1 \otimes P + P \otimes 1 \tag{1.6}$$

The primitives of a Hopf algebra form a Lie algebra under the commutator product

$$[P_1, P_2] = P_1 P_2 - P_2 P_1 \tag{1.7}$$

which is denoted Prim($H$). For any Hopf algebra there is strong interest ina description of its Lie algebra of primitives. For instance because of the Milnor - Moore theorem, [36], that says that a graded connected cocommutative Hopf algebra over a field of characteristic zero is isomorphic to to the universal enveloping algebra of its Lie algebra of primitives. Also, far from unrelated, let $Q(H) = I(H)/I(H)^2$ be the module of indecomposables of a graded Hopf algebra $H$. Here $I(H)$ is the augmentation ideal of $H$. Then there is an induced duality between $Q(H)$ and Prim($H^*$), and there is the (classical) Leray theorem that says that for a connected commutative graded Hopf algebra $H$ over a characteristic zero field any section of $I(H) \longrightarrow Q(H)$ induces an isomorphism of the free commutative algebra over $Q(A)$ to $H$. This last theorem now has been considerably generalized to the setting of operads, see [38], and the references quoted there.

The first main topic that is treated in some detail (but without proofs) in this survey is an explicit and algorithmic description of a basis over the integers of Prim(*NSymm*).

A *divided power sequence* in a Hopf algebra $H$ is a sequence of elements

$$d = (d(0) = 1, d(1), d(2), \cdots) \tag{1.8}$$

such that for all $n$

$$\mu_H(d(n)) = \sum_{i+j=n} d(i) \otimes d(j) \quad i, j \in \{1, 2, 3, \cdots\} \tag{1.9}$$

Note that $d(1)$ is a primitive. Is is sometimes useful to write a DPS (divided power sequence) as a power series in a counting variable $t$:

$$d(t) = 1 + d(1)t + d(2)t^2 + d(3)t^3 + \cdots \tag{1.10}$$

That makes it easier to talk about the inverse of a DPS (inverse power series), the product of two DPS's (multiplication of power series) and shifted DPS's: $d(t) \mapsto d(t^n)$, all operations that give new DPS's from old ones. When written in the form (1.10) a DPS is often called a *curve*.

It turns out that each primitive of Prim(*NSymm*) can be extended to a divided power sequence. This is important because it implies that as a coalgebra *NSymm* is the cocommutative cofree graded coalgebra over the module Prim(*NSymm*).

Now let *QSymm* be the graded dual Hopf algebra (over the integers) of *NSymm*. For an explicit description of *QSymm*, the Hopf algebra of quasisymmetric functions, see below in section 2.
    A most important question concerning *QSymm* is whether it is free polynomial as a commutative algebra. This has been an important issue since 1972, since it is crucial for the development of certain parts of the theory of noncommutative formal groups, [5, 6, 17]. The matter was finally settled in 1999, [21], in the positive sense that it is indeed free. A second proof follows from the cofreeness of *NSymm*. However, both these proofs fail to produce explicit generators. This has now also been taken care of, [23], and is the second main topic that will be discussed in some detail below.

One most interesting and important aspect of the structure of *Symm* is the presence of two families of Hopf algebra morphisms that are called Frobenius and Verschiebung morphisms. They satisfy a large number of beautiful relations. The third main topic of this survey is to what extent these can be lifted to *NSymm*, respectively, extended to *QSymm*. There are both positive and negative results. However, the matter has not yet been quite completely settled.

This paper is an expanded write-up of two talks that I gave on the subject: in Krasnoyarsk in August 2002 at the occasion of the International Conference "Algebra and its applications" in honour of the 70-th anniversary of V P Shunkov and the 65-th anniversary of V M Busarkin, and at the Z. Borewicz memorial conference in Skt Petersburg in September 2002.

**2. The Hopf algebra** *QSymm* **of quasisymmetric functions.**

Above, in the introduction, the graded Hopf algebra *NSymm* of noncommutative symmetric functions was defined. The grading is defined by

$$\text{wt}(Z_n) = n \tag{2.1}$$

and, more generally, if $\alpha = [a_1, a_2, \cdots, a_m]$ is a nonempty word over the positive integers $\mathbf{N} = \{1, 2, \cdots\}$, let $Z_\alpha$ be the noncomutative monomial

$$Z_\alpha = Z_{a_1} Z_{a_2} \cdots Z_{a_m} \tag{2.2}$$

then

$$\text{wt}(Z_\alpha) = \text{wt}(\alpha) = a_1 + \cdots + a_m \tag{2.3}$$

Let $Z_{[]} = 1$, where [ ] is the empty word, then the $Z_\alpha$, $\alpha \in \mathbf{N}^*$, the monoid of words over $\mathbf{N}$ form a basis of *NSymm* (as a graded Abelian group). The empty word, and also $Z_{[]} = 1$, has weight zero.

As a free Abelian graded group *QSymm*, the graded dual of *NSymm* can be taken to be the free Abelian group with as basis $\mathbf{N}^*$, the words over the set of natural numbers. The duality is then

$$<Z_\alpha, \beta> = \delta_\alpha^\beta \tag{2.4}$$

The duality induced comultiplication is easy to describe. It is '*cut*':

$$[a_1, a_2, \cdots, a_m] \mapsto \sum_{i=0}^{m} [a_1, \cdots, a_i] \otimes [a_{i+1}, \cdots, a_m] \tag{2.5}$$

where of course $[a_1, \cdots, a_i] = [\ ] = 1$ if $i = 0$ and $[a_{i+1}, \cdots a_m] = [\ ] = 1$ if $i = m$. The duality induced multiplication is more difficult to describe. It is the socalled '*overlapping shuffle multiplication*' which can be described as follows.

Let $\alpha = [a_1, a_2, \cdots, a_m]$ and $\beta = [b_1, b_2, \cdots, b_n]$ be two compositions or words. Take a 'sofar empty' word with $n + m - r$ slots where $r$ is an integer between 0 and $\min\{m, n\}$, $0 \le r \le \min\{m, n\}$. Choose $m$ of the available $n + m - r$ slots and place in it the natural numbers from $\alpha$ in their original order; choose $r$ of the now filled places; together with the remaining $n + m - r - m = n - r$ places these form $n$ slots; in these place the entries from $\beta$ in their orginal order; finally, for those slots which have two entries, add them. The product of the two words $\alpha$ and $\beta$ is the sum (with multiplicities) of all words that can be so obtained. So, for instance,

$$[a,b] \times_{osh} [c,d] = [a,b,c,d] + [a,c,b,d] + [a,c,d,b] + [c,a,b,d] + [c,a,d,b] + [c,d,a,b] + $$
$$+ [a+c,b,d] + [a+c,d,b] + [c,a+d,b] + [a,b+c,d] + [a,c,b+d] + \tag{2.6}$$
$$+ [c,a,b+d] + [a+c,b+d]$$

and $[1] \times_{osh} [1] \times_{osh} [1] = 6[1,1,1] + 3[1,2] + 3[2,1] + [3]$.

There is a concrete realization of *QSymm* much like the standard realization of *Symm* as the ring of symmetric functions in infinitely many indeterminates $x_1, x_2, \cdots$. See [34], Chapter 1 for some detail on how to work with infinitely many indeterminates in this context.

Let $X$ be a finite or infinite set (of commuting variables) and consider the ring of polynomials, $R[X]$, and the ring of power series, $R[[X]]$, over a commutative ring $R$ with unit element in the commuting variables from *X*. A polynomial or power series $f(X) \in R[[X]]$ is called *symmetric* if for any two finite sequences of indeterminates $x_1, x_2, \cdots, x_n$ and $y_1, y_2, \cdots, y_n$ from *X* and any sequence of exponents $i_1, i_2, \cdots, i_n \in \mathbf{N}$, the coefficients in $f(X)$ of $x_1^{i_1} x_2^{i_2} \cdots x_n^{i_n}$ and $y_1^{i_1} y_2^{i_2} \cdots y_n^{i_n}$ are the same.

The *quasi-symmetric formal power series* are a generalization introduced by Gessel, [15], in connection with the combinatorics of plane partitions. This time one takes a *totally ordered* set of indeterminates, e.g. $V = \{v_1, v_2, \cdots\}$, with the ordering that of the natural numbers, and the condition is that the coefficients of $x_1^{i_1} x_2^{i_2} \cdots x_n^{i_n}$ and $y_1^{i_1} y_2^{i_2} \cdots y_n^{i_n}$ are equal for all totally ordered sets of indeterminates $x_1 < x_2 < \cdots < x_n$ and $y_1 < y_2 < \cdots < y_n$. Thus, for example,

$$x_1 x_2^2 + x_2 x_3^2 + x_1 x_3^2 \tag{2.7}$$

is a quasi-symmetric polynomial in three variables that is not symmetric.

Products and sums of quasi-symmetric polynomials and power series are again quasi-symmetric (obviously), and thus one has, for example, the ring of quasi-symmetric power series $QSymm^\wedge$ in countably many commuting variables over the integers and its subring

$$QSymm \tag{2.8}$$

of quasi-symmetric polynomials in finite of countably many indeterminates, which are the quasi-symmetric power series of bounded degree. The notation is justified. The quasisymmetric functions in $\{x_1, x_2, \cdots\}$ in this sense are a concrete realization of the quasisymmetric functions as introduced above as the graded dual of *NSymm*.

In detail, given a word $\alpha = [a_1, a_2, \cdots, a_m]$ over $\mathbf{N}$, also called a *composition* in this context, consider the *quasi-monomial function*

$$M_\alpha = \sum_{i_1 < \cdots < i_m} x_{i_1}^{a_1} x_{i_2}^{a_2} \cdots x_{i_m}^{a_m} \tag{2.9}$$

defined by $\alpha$. It is now an easy exercise to verify that as power series in the $\{x_1, x_2, \cdots\}$ the $M_\alpha$ satisfy

$$M_\alpha M_\beta = M_{\alpha \times_{osh} \beta} \tag{2.10}$$

where $\times_{osh}$ is the overlapping shuffle product of words just defined above (2.6).

## 3. *NSymm* and *QSymm* over a field of characteristic zero.

Let *K* be a field of characteristic zero, in particular the field of rational numbers $\mathbf{Q}$. The Hopf algebras of noncommutative symmetric functions and quasisymmetric functions over *K* are denoted

$$NSymm_K = NSymm \otimes_{\mathbf{Z}} K, \quad QSymm_K = QSymm \otimes_{\mathbf{Z}} K \tag{3.1}$$

As remarked in the introduction, working over a field of characteristic zero tends to simplify things considerably. Not that then everything becomes clear and easy and straightforward. Very far from it; witness the many papers quoted in the introduction. However, it is certainly true, that for the three groups of questions which form the main topic of this paper: primitives of noncommutative symmetric functions, freeness of the algebra of quasisymmetric functions, existence of Frobenius and Verschiebung morphisms, things either reduce to known things, or become fairly straightforward. At least up to fairly explicitly given isomorphisms. And in that connection it is well to reflect that knowing something well up to an isomorphism can be not all the same as really controlling things.

To start with, consider another Hopf algebra over the integers, the *Lie Hopf algebra*

$$\mathcal{U} = \mathbf{Z}\langle U_1, U_2, \cdots \rangle, \quad \mu(U_n) = 1 \otimes U_n + U_n \otimes 1, \quad \varepsilon(U_n) = 0, \quad n = 1, 2, \cdots \quad (3.2)$$

That is, as an algebra $\mathcal{U}$ is the free algebra over the integers in the noncommuting variables $U_n$ and the coalgebra structure is given by requiring that the two right hand formulas of (3.2) are algebra homomorphisms.

The primitives of $\mathcal{U}$ are called *Lie polynomials* and they form the free Lie algebra over $\mathbf{Z}$ in the alphabet $\{U_1, U_2, \cdots\}$. Also the universal enveloping algebra of $\mathrm{Prim}(\mathcal{U})$ is $\mathcal{U}$.

Now, over the rationals (and hence over any field of characteristic zero or ring containing $\mathbf{Q}$) $NSymm$ [1] and $\mathcal{U}$ are isomorphic. One particularly beautiful isomorphism $NSymm_\mathbf{Q} \xrightarrow{\varphi} \mathcal{U}_\mathbf{Q}$ is given by setting

$$\begin{aligned} 1 + Z_1 t + Z_2 t^2 + Z_3 t^3 + \cdots &= \exp(U_1 t + U_2 t^2 + U_3 t^3 + \cdots) \\ &= \sum_{i=0}^{\infty} (U_1 t + U_2 t^2 + U_3 t^3 + \cdots)^i \end{aligned} \quad (3.3)$$

where $t$ is a (counting) variable commuting with everything in sight. Equating equal powers of $t$ on the left and right hand sides of (3.3) gives

$$Z_n = \sum_{\substack{i_1 + \cdots + i_k = n \\ i_j \in \mathbf{N}}} \frac{U_{i_1} U_{i_2} \cdots U_{i_k}}{k!} \quad (3.4)$$

For a proof that the algebra isomorphism defined by setting $\varphi(Z_n)$ equal to the right hand side of (3.4) defines an isomorphism of Hopf algebras see [18].

Thus, up to the isomorphism $\varphi$, the matter of describing $\mathrm{Prim}(NSymm)$, comes down to describing $\mathrm{Prim}(\mathcal{U})$, or, equivalently writing down some explicit bases for the free Lie algebra generated by a countable set of indeterminates $\{U_1, U_2, \cdots\}$. The matter of constructing bases for free Lie algebras has had a lot af attention. There are many more or less different bases such as Hall bases, Shirshov bases, Lyndon bases, see [3, 39, 41, 42, 47].

Here, partly because most of the concepts will be needed below anyway, is a description of the *Lyndon basis* (also callled *Chen - Fox - Lyndon basis*).

Consider the free monoid $\mathbf{N}^*$ of words on the alphabet $\mathbf{N}$ (or any other totally ordered alphabet for that matter). The *lexicographic order* (also called dictionary order, or alphabetical order) on $\mathbf{N}^*$ is defined as follows. If $\alpha = [a_1, a_2, \cdots, a_m]$ and $\beta = [b_1, b_2, \cdots, b_n]$ are two words of length $m$ and $n$ respectively, $\mathrm{lg}(\alpha) = m$, $\mathrm{lg}(\beta) = n$,

$$\alpha >_{lex} \beta \iff \begin{cases} \exists k \leq \min\{m, n\} \text{ such that } a_1 = b_1, \cdots, a_{k-1} = b_{k-1} \text{ and } a_k > b_k \\ \text{or } \mathrm{lg}(\alpha) = m > \mathrm{lg}(\beta) = n \text{ and } a_1 = b_1, \cdots, a_n = b_n \end{cases} \quad (3.5)$$

The empty word is smaller than any other word. This defines a total order. Of course, if one accepts the dictum that anything is larger than nothing, the second clause of (3.5) is superfluous.

The *proper tails* (*suffixes*) of the word $\alpha = [a_1, a_2, \cdots, a_m]$ are the words $[a_i, a_{i+1}, \cdots a_m]$, $i = 2, 3, \cdots, m$. Words of length 1 or 0 have no proper tails. The *prefix* corresponding to a tail $\alpha'' = [a_i, a_{i+1}, \cdots a_m]$ is $\alpha' = [a_1, \cdots, a_{i-1}]$ so that $\alpha = \alpha' * \alpha''$ where $*$

---

[1] The Hopf algebra *NSymm* is sometimes called the Leibniz Hopf algebra.

denotes *concatenation* of words.

A word is *Lyndon* iff it is lexicographically smaller than each of its proper tails. For instance [4], [1,3,2], [1,2,1,3] are Lyndon and [1,2,1] and [2,1,3] are not Lyndon.

For each Lyndon word $\alpha$ of length $>1$ consider the lexicographically smallest proper tail $\alpha''$ of $\alpha$. Let $\alpha'$ be the corresponding prefix to $\alpha''$. Then $\alpha'$ and $\alpha''$ are both Lyndon and $\alpha = \alpha' * \alpha''$ is called the *canonical factorization* of $\alpha$.

A basis of the free Lie algebra on $\{U_1, U_2, \cdots\}$, i.e. a basis of $\mathrm{Prim}(\mathcal{U}) \subset \mathcal{U}$, is now obtained as follows. For each word $\alpha = [a_1, a_2, \cdots, a_m]$ let $U_\alpha = U_{a_1} U_{a_2} \cdots U_{a_m}$ be the corresponding monomial. Now, by recursion in length, define for a word of length 1

$$Q_{[i]} = U_i \qquad (3.6)$$

and for $\alpha$ Lyndon and of length $\lg(\alpha) \geq 2$ let $\alpha = \alpha' * \alpha''$ be its canonical factorization and set

$$Q_\alpha = [Q_{\alpha'}, Q_{\alpha''}] \qquad (3.7)$$

then the $\{Q_\alpha : \alpha \text{ Lyndon}\}$ form a basis of $\mathrm{Prim}(\mathcal{U}) \subset \mathcal{U}$. For a proof see e.g. [39], p. 105ff.

The next topic to be taken up is the matter of the freeness of *QSymm* over the rationals. The graded dual of $\mathcal{U}$ is the socalled *shuffle algebra*. As a free module over $\mathbf{Z}$ it has the words over $\mathbf{N}$ as a basis and the product is the *shuffle product* which is like the overlapping shuffle product except that the overlap terms, i.e. those which involve additions of entries are left out. Thus for example

$$[a,b] \times_{sh} [c,d] = [a,b,c,d] + [a,c,b,d] + [a,c,d,b] + [c,a,b,d] + [c,a,d,b] + [c,d,a,b]$$

(compare (2.6) above).

It is well known that the shuffle algebra is free polynomial with as generators (for example) the Lyndon words. See, for example, [39]. p. 111 for a proof. Thus via the isomorphism $\varphi$, or rather its graded dual, it follows that $QSymm_\mathbf{Q}$ is a free commutative algebra. But the description of the generators is rather involved and they do not look very nice. Actually the situation is rather better and a modification of the proof of the freeness of the shuffle algebra (using a different ordering on words) gives that in fact $QSymm_\mathbf{Q}$ is commutative free polynomial on the Lyndon words. The ordering to be used is the wll-ordering. The acronym stands for weight first, than length, than lexicographic. See [20] for details.

The third main topic of this survey is the existence of Frobenius and Verschiebung type Hopf algebra endomorphisms of *NSymm* and *QSymm* which lift, respectively extend, those on *Symm*. Again, over the rationals, this is a relatively straightforward matter. Though there are some unanswered questions.

Recall the situation for *Symm*, see [17, 24] for more details.

On *Symm* there are two families of Hopf algebra endomorphims, called Frobenius and Verschiebung morphisms, denoted $\mathbf{f}_n$, $\mathbf{v}_n$, $n \in \mathbf{N}$, which among others have the following beautiful properties:

(i) $\mathbf{f}_1 = \mathbf{v}_1 = id$

(ii) $\mathbf{f}_n$ is homogeneous of degree $n$, i.e. $\mathbf{f}_n(Symm_k) \subset Symm_{nk}$

Here, for any graded Hopf algebra, $H$, $H_n$ is the homogeneous part of of weight $n$ of $H$.

(iii) $\mathbf{v}_n$ is homogenous of degree $n^{-1}$, i.e. $\mathbf{v}_n(Symm_k) \subset Sym_{n^{-1}k}$ if $n$ divides $k$, and $\mathbf{v}_n(Symm_k) = 0$ if $n$ does not divide $k$.

(iv) $\mathbf{f}_n \mathbf{f}_m = \mathbf{f}_{nm}$ for all $n, m \in \mathbf{N}$

(v)   $\mathbf{v}_n \mathbf{v}_m = \mathbf{v}_{nm}$  for all  $n, m \in \mathbf{N}$
(vi)  $\mathbf{f}_n \mathbf{v}_m = \mathbf{v}_m \mathbf{f}_n$  provided  $n$  and  $m$  are relatively prime,  $\gcd(m, n) = 1$
(vii) $\mathbf{v}_n \mathbf{f}_n = \mathbf{n}$, where  $\mathbf{n}$  is the $n$-fold convolution of the identiy.

Now there is the natural projection

$$NSymm \longrightarrow Symm, \quad Z_n \mapsto h_n \tag{3.8}$$

and the natural (graded dual) inclusion

$$Symm \subset QSymm \tag{3.9}$$

obtained by regarding a symmetric function as a special kind of quasisymmetric function. The question is whether there are lifts, respectively extensions, on *NSymm*, respectively *QSymm*, which also have the properties (i) - (vii).

Retaining property (vii) can be ruled out immediately for trivial reasons. The simple fact is that **n** on either *Qsymm* or *NSymm* simply is not a Hopf algebra endomorphism. So it is natural to concentrate on the other six properties. And then the answer over the rationals is yes. But, as will be stated below, the answer over the integers is no. But there are interesting substitutes.

Let

$$p_n = x_1^n + x_2^n + x_3^n \cdots \tag{3.10}$$

denote the power sums in *Symm*. They are related to the complete symmetric functions by the recursion relation

$$nh_n = p_n + p_{n-1}h_1 + p_{n-2}h_2 + \cdots p_1 h_{n-1} \tag{3.11}$$

The *Frobenius* and *Verschiebung* morphisms on *Symm* are characterized by

$$\mathbf{f}_n p_k = p_{nk}, \quad \mathbf{v}_n p_k = \begin{cases} np_{k/n} & \text{if } n \text{ divides } k \\ 0 & \text{if } n \text{ does not divide } k \end{cases} \tag{3.12}$$

On the polynomial generators  $h_n$  this characterization of  $\mathbf{v}_n$  works out as

$$\mathbf{v}_n h_k = \begin{cases} h_{k/n} & \text{if } n \text{ divides } k \\ 0 & \text{otherwise} \end{cases} \tag{3.13}$$

Define the (noncommutative) *Newton primitives* in *NSymm* by

$$P_n(Z) = \sum_{r_1 + \cdots r_k = n} (-1)^{k+1} r_k Z_{r_1} Z_{r_2} \cdots Z_{r_k}, \quad r_i \in \mathbf{N} = \{1, 2, \cdots\} \tag{3.14}$$

or, equivalently, by the recursion relation

$$nZ_n = P_n(Z) + Z_1 P_{n-1}(Z) + Z_2 P_{n-2}(Z) + \cdots + Z_{n-1} P_1(Z) \tag{3.15}$$

Note that under the projection  $Z_n \mapsto h_n$  by (3.15) and (3.11)  $P_n(Z)$  goes to  $p_n$. It is easily proved by induction, using (3.15), or directly from (3.14), that the  $P_n(Z)$  are primitives of *NSymm*, and it is also easy to see from (3.15) that over the rationals *NSymm* is the free associative algebra generated by the $P_n(Z)$. Thus over the rationals the Lie algebra of primitives of *NSymm* is simply the free Lie algebra generated by the $P_n(Z)$, giving a second description of Prim( $NSymm_\mathbf{Q}$ ).

There are obvious candidate lifts of the $\mathbf{v}_n$ on *Symm* to Hopf algebra endomorphisms on *NSymm.*, viz

$$\mathbf{v}_n(Z_k) = \begin{cases} Z_{k/n} & \text{if } k \text{ is divisible by } n \\ 0 & \text{otherwise} \end{cases} \tag{3.16}$$

By (3.14) or (3.15) this implies

$$\mathbf{v}_n(P_k) = \begin{cases} nP_{k/n} & \text{if } n \text{ divides } k \\ 0 & \text{otherwise} \end{cases} \tag{3.17}$$

Now on $NSymm_\mathbf{Q}$ define the Frobenius morphisms as the algebra morphisms given by

$$\mathbf{f}_n(P_k(Z)) = P_{nk}(Z) \tag{3.18}$$

It is now easily checked that the $\mathbf{v}_n$ and $\mathbf{f}_n$ as defined by (3.16) and (3.18) are Hopf algebra endomorphisms of $NSymm_\mathbf{Q}$ , that they satisfy (the analogues on $NSymm_\mathbf{Q}$ of) properties (i)-(vi) and that they descend to the usual Frobenius and Verschiebung morphisms on *Symm.*

A priori, the $\mathbf{f}_n$ as defined by (3.18) are only defined over the rationals and indeed nontrivial denominators show up almost immediately. For instance

$$\begin{aligned} \mathbf{f}_2(Z_1) &= 2Z_2 - Z_1^2 \\ \mathbf{f}_2(Z_2) &= 2Z_4 - \tfrac{3}{2}Z_1Z_3 - \tfrac{1}{2}Z_3Z_1 + Z_2^2 + \tfrac{1}{2}Z_1Z_2Z_1 + \tfrac{1}{2}Z_1^2Z_2 \end{aligned} \tag{3.19}$$

On *Symm* a certain amount of coefficient magic sees to it that all coefficients become integral. But of course over *Symm* there are much better definitions of the Frobenius morphisms that immediately show that they are defined over the integers, see [24] or [17], §17.

As we shall see later, over the integers there are even no algebra endomorphisms $\mathbf{f}_n$ of *NSymm* that lift the $\mathbf{f}_n$ on *Symm* such that together with the $\mathbf{v}_n$ as defined by (3.16) they satisfy (i)-(vi).

Note there is nothing unique about this solution (3.18) of the Frobenius-Verschiebung lifting problem over the rationals. For instance one could work instead with the seond set of Newton primitives defined by

$$P'_n(Z) = \sum_{r_1 + \cdots + r_k = n} (-1)^{k+1} r_1 Z_{r_1} Z_{r_2} \cdots Z_{r_k}, \quad r_i \in \mathbf{N} = \{1,2,\cdots\} \tag{3.20}$$

and satisfying the recursion relation

$$nZ_n = P'_n(Z) + P'_{n-1}(Z)Z_1 + P'_{n-2}(Z)Z_2 + \cdots + P'_1(Z)Z_{n-1} \tag{3.20}$$

**4. The primitives of *NSymm*.**
Above, some primitives of *NSymm* were already written down and they generate a free graded Lie algebra contained in Prim(*NSymm*). Denote this Lie algebra by FrLie(*P*) and its homogeneous part of weight $n$ by FrLie(*P*)$_n$. The Lie algebra Prim(*NSymm*) is also graded of course. Let $\mathrm{Prim}(NSymm)_n \subset NSymm_n$ be the homogeneous part of degree $n$ of Prim($NSymm$). Both $\mathrm{Prim}(NSymm)_n$ and $\mathrm{FrLie}(P)_n$ are free Abelian groups of rank $\beta_n$, the number of weight $n$ Lyndon words [2]. The index of $\mathrm{FrLie}(P)_n \subset \mathrm{Prim}(NSymm)_n$ as a function of $n$ measures how large $\mathrm{FrLie}(P)$ is in $\mathrm{Prim}(NSymm)$. As it turns out $\mathrm{FrLie}(P)$ is only a tiny part. Indeed, the value of the index alluded to is

$$\text{Index of FrLie}(P)_n \text{ in Prim}(NSymm)_n = \prod_{\alpha \in LYN, \text{wt}(\alpha)=n} \frac{k(\alpha)}{g(\alpha)} \quad (4.1)$$

where for a word $\alpha = [a_1, a_2, \cdots, a_m]$ over the natural numbers $g(\alpha)$ is the gcd (greatest common divisor) of its entries $a_1, a_2, \cdots, a_m$ and $k(\alpha)$ is the product of its entries. Thus the values of (4.1) for the first six $n$ are 1, 1, 2, 6, 576, 69120.

Thus taking iterated commutators of the known Newton primitives is not nearly good enough. One can see this coming very quickly. Indeed $[P_1, P_2] = 2(Z_1 Z_2 - Z_2 Z_1)$. It also follows that $\mathrm{Prim}(NSymm)$ is not a free Lie algebra over the integers. Rather it tries to be something like a free divided power Lie algebra (though I do not know what such a thing might be).

Instead of taking commutators of primitives it turns out to be a good idea to work with whole DPS's (divided power sequences, see (1.8)). Accordingly, the next thing to be described are techniques for producing new divided power sequences from known ones. There are two more techniques for this (besides the ones mentioned in the introduction, which do not suffice) coming from two socalled isobaric decomposition theorems.

For the first *isobaric decomposition theorem* consider the Hopf algebra

$$2NSymm = \mathbf{Z}\langle X_1, Y_1, X_2, Y_2, \cdots\rangle, \ \mu(X_n) = \sum_{i+j=n} X_i \otimes X_j, \ \mu(Y_n) = \sum_{i+j=n} Y_i \otimes Y_j \quad (4.2)$$

and the two natural curves

$$X(s) = 1 + X_1 s + X_2 s^2 + \cdots, \quad Y(t) = 1 + Y_1 t + Y_2 t^2 + \cdots \quad (4.3)$$

and consider the commutator product

$$X(s)^{-1} Y(t)^{-1} X(s) Y(t) \quad (4.4)$$

On the set of pairs of nonnegative integers consider the ordering

$$(u, v) <_{wl} (u', v') \iff u + v < u' + v' \text{ or } (u + v = u' + v' \text{ and } u < u') \quad (4.5)$$

(Here the index $wl$ on $<_{wl}$ is supposed to be a mnemonic for weight first, then lexicographic.)

4.6. *Theorem* (first bi-isobaric decomposition theorem, Shay [40], Ditters). There are 'higher commutators' (or perhaps better 'corrected commutators')

$$L_{u,v}(X, Y) \in \mathbf{Z}\langle X, Y\rangle, \ (u,v) \in \mathbf{N} \times \mathbf{N} \quad (4.7)$$

such that

---

[2] The numbers $\beta_n$ are given by the identity $(1-t)^{-1}(1-2t) = \prod_{n=1}^{\infty}(1-t^n)^{\beta_n}$ which goes back to Witt, [48].

$$X(s)^{-1}Y(t)^{-1}X(s)Y(t) = \overrightarrow{\prod_{\gcd(a,b)=1}}(1+ L_{a,b}(X,Y)s^a t^b + L_{2a,2b}(X,Y)s^{2a}t^{2b} +\cdots) \quad (4.8)$$

where the product is an ordered product for the ordering $<_{wl}$ just introduced, (4.5). Moreover

(i) $L_{u,v}(X,Y) = [X_u, Y_v] + $ (terms of length $\geq 3$) \hfill (4.9)

(ii) $L_{u,v}(X,Y)$ is homogeneous of weight $u$ in $X$ and of weight $v$ in $Y$. \hfill (4.10)

(iii) For $\gcd(a,b)=1$, $1+ L_{a,b}(X,Y)s^a t^b + L_{2a,2b}(X,Y)s^{2a}t^{2b} +\cdots$ is a 2-curve.

Here a 2-curve is a two dimensional version of a curve. A power series in two variables with constant term 1

$$d(s,t) = \sum_{i,j} d(i,j) s^i t^j \quad (4.11)$$

is a 2-curve iff

$$\mu(d(m,n)) = \sum_{\substack{m_1+m_2=m \\ n_1+n_2=n}} d(m_1,n_1) \otimes d(m_2,n_2) \quad (4.12)$$

This is not at all difficult to prove. The only thing needed is to observe that pure powers of $s$ or $t$ do not occur on the LHS of (4.8) and that each pair of nonnegative integers $(u,v)$ occurs just once in one of the factors on the right 0f (4.8). That gives the decomposition. The fact that the factors are two curves then follows easily with induction from the observation that the LHS of (4.8) is a 2-curve.

Also (4.8) implies an explicit recursion formula for the $L_{u,v}(X,Y)$.

4.13. *Theorem* (Second bi-isobaric decomposition theorem, Hazewinkel [22]). There are unique homogeneous noncommutative polynomials $N_{u,v}(Z) \in NSymm$ such that

$$Z(s)^{-1}Z(t)^{-1}Z(s+t) = \overrightarrow{\prod_{\substack{a,b\in \mathbf{N} \\ \gcd(a,b)=1}}}(1+ N_{a,b}(Z)s^a t^b + N_{2a,2b}(Z)s^{2a}t^{2b} +\cdots). \quad (4.14)$$

Moreover

(i) $N_{u,v}(Z) = \binom{u+v}{u} Z_{u+v} + $ (terms of length $\geq 2$) \hfill (4.15)

(ii) $N_{u,v}(Z)$ is homogeneous of weight $u+v$ \hfill (4.16)

(iii) For each $a,b \in \mathbf{N}^2$, $\gcd(a,b)=1$,

$$1+ N_{a,b}(Z)s^a t^b + N_{2a,2b}(Z)s^{2a}t^{2b} +\cdots \quad (4.17)$$

is a 2-curve.

(iv) For each $n\geq 2$, $N_{1,n-1}(Z) = P_n(Z)$ \hfill (4.18)

Again, the decomposition is not at all difficult to prove and the final observation results directly from the recursion formula implied by (4.14) compared to the recursion formula (3.15) for the $P_n(Z)$.

There is now a sufficiency of tools to describe a basis of Prim(*NSymm*) and more. The essential fact is that if $d_1, d_2$ are two divided power sequences in some Hopf algebra $H$, than so is

$$D_{a,b}(d_1, d_2) = (1, L_{a,b}(d_1, d_2), L_{2a,2b}(d_1, d_2), \cdots) \tag{4.19}$$

where, as the notation suggests, $L_{u,v}(d_1, d_2)$ is obtained from the $L_{u,v}(X,Y)$ of theorem 4.6 by substituting $d_1(k)$ for $X_k$ and $d_2(l)$ for $Y_l$. This follows immediately from the fact that for $\gcd(a,b) = 1$

$$1 + L_{a,b}(X,Y)s^a t^b + L_{2a,2b}(X,Y)s^{2a} t^{2b} + \cdots$$

is a 2-curve. Similarly, if $d$ is a curve in any Hopf algebra $H$

$$N_{a,b}(d) = 1 + N_{a,b}(d)t + N_{2a,2b}(d)t^2 + \cdots \tag{4.20}$$

is another curve.³

Let *LYN* denote the set of Lyndon words over the natural numbers $\mathbf{N} = \{1,2,3,\cdots\}$. Then to each $\alpha = [a_1, \cdots, a_m] \in LYN$ there are associated three things

(i)　A number $g(\alpha) = \gcd\{a_1, \cdots, a_m\}$
(ii)　A divided power sequence $d_\alpha$
(iii)　A primitive $P_\alpha$

The items (ii) and (iii) are defined recursively as follows. If $\lg(\alpha) = m$ is 1, $d_{[n]} = (1, Z_1, Z_2, \cdots)$ and $P_{[n]} = P_n(d_{[n]}) = P_n(Z)$. If $\lg(\alpha) \geq 2$, let $\alpha = \alpha' * \alpha''$ be the canonical factorization of $\alpha$ (see just above (3.6)). Then

$$\begin{aligned}d_\alpha &= (1, L_{g(\alpha')/g(\alpha), g(\alpha'')/g(\alpha)}(d_{\alpha'}, d_{\alpha''}), L_{2g(\alpha')/g(\alpha), 2g(\alpha'')/g(\alpha)}(d_{\alpha'}, d_{\alpha''}), \cdots) \\ &= D_{g(\alpha')/g(\alpha), g(\alpha'')/g(\alpha)}(d_{\alpha'}, d_{\alpha''})\end{aligned} \tag{4.21}$$

and

$$P_\alpha = P_{g(\alpha)}(d_\alpha) \tag{4.22}$$

Note that the divided power sequences associated to $[a_1, \cdots, a_m]$ and $[ra_1, \cdots, ra_m]$, $r \in \mathbf{N}$ are the same.

　　4.23. *Theorem.* The $P_\alpha$, $\alpha \in LYN$ form a basis over the integers of Prim(*NSymm*). Each of the $P_\alpha$ is the first term of a DPS.

The second property of theorem 4.23 guarantees that *NSymm* is the cofree cocommutative graded coalgebra over the graded module Prim(*NSymm*), see the appendix of [22] for a proof of that. It follows that the graded dual *QSymm* is a free commutative algebra, and implicitly specifies a set of generators for *QSymm*. However, this does not give a convenient description of such a set of generators.

---

　　³ The operations on curves defined by (4.19) and (4.20) are functorial. Thus assigning to a Hopf algebra $H$ its group of curves gives a group valued functor with many operations on it, viz 2-ary operations $L_{a,b}$ for every pair $a,b$ with $\gcd(a,b)=1$, and also unary operations $N_{a,b}$ for any such pair. In addition there are the operations which come from the Verschiebung Hopf algebra endomorphisms on *NSymm* (see section 6) which commute with the $L_{a,b}$ and $N_{a,b}$. And then there are also Frobenius type operations (see also section 6). It is almost totally uninvestigated whether some subset of all this structure can be used for classification purposes, say, of noncommutative formal groups.

The original proof of theorem 4.23, [22], is rather long and intricate. Fortunately there is now a much shorter proof, which will be discussed in the next section.

The basis $P_\alpha$, $\alpha \in LYN$ of Prim(*NSymm*) has a number of nice properties, particularly with respect to the Verschiebung endomorphisms $\mathbf{v}_n$ of *NSymm*, see (3.16).

Consider the following ordering on words: $\alpha <_{wll} \beta$ if and only if ($\text{wt}(\alpha) < \text{wt}(\beta)$ or ($\text{wt}(\alpha) = \text{wt}(\beta)$ and $\lg(\alpha) < \lg(\beta)$) or ($\text{wt}(\alpha) = \text{wt}(\beta)$ and $\lg(\alpha) = \lg(\beta)$) and $\alpha <_{lex} \beta$)). The acronym 'wll' stands for 'weight first, than length, then lexicographic'.

4.24. *Theorem.* For each Lyndon word $\alpha = [a_1, \cdots, a_m]$

(i) $P_\alpha = g(\alpha) Z_\alpha +$ (wll larger terms) \hfill (4.25)

(ii) $\mathbf{v}_n(P_\alpha) = \begin{cases} nP_{[n^{-1}a_1, \cdots, n^{-1}a_m]} & \text{if } n \text{ divides } g(\alpha) \\ 0 & \text{otherwise} \end{cases}$ \hfill (4.26)

The first property of theorem 4.24 results directly from the construction. The second property comes from the fact that the curves $\sum L_{ra,rb}(X,Y) s^{ra} t^{rb}$ and $\sum N_{ra,rb}(Z) s^{ra} t^{rb}$ of the two isobaric decomposition theorems 4.6 and 4.13 are particularly nice curves in that

$$\mathbf{v}_n(L_{u,v}(X,Y)) = \begin{cases} L_{n^{-1}u, n^{-1}v}(X,Y) & \text{if } n \text{ divides } u,v \\ 0 & \text{otherwise} \end{cases} \quad (4.27)$$

$$\mathbf{v}_n(N_{u,v}(Z)) = \begin{cases} N_{n^{-1}u, n^{-1}v}(Z) & \text{if } n \text{ divides } u,v \\ 0 & \text{otherwise} \end{cases} \quad (4.28)$$

where $\mathbf{v}_n$ is defined on *2NSymm* just like (3.16) for both the *X*'s and the *Y*'s.

## 5. Free polynomial generators for *QSymm* over the integers.

As in section 2 above, consider *Symm* and *QSymm* as symmetric functions in an infinity of ideterminates

$$Symm = \mathbf{Z}[e_1, e_2, \cdots] \subset QSymm \subset \mathbf{Z}[x_1, x_2, \cdots] \quad (5.1)$$

where the $e_i$ are the elementary symmetric functions in the $x_j$. There is a well known $\lambda$-ring structure on $\mathbf{Z}[x_1, x_2, \cdots]$ given by

$$\lambda_i(x_j) = \begin{cases} x_j & \text{if } i = 1 \\ 0 & \text{if } i \geq 2 \end{cases}, \quad j = 1, 2, \cdots \quad (5.2)$$

For information on $\lambda$-rings see [27]. The asscociated *Adams operators*, determined by the formula

$$t \frac{d}{dt} \log \lambda_t(a) = \sum_{n=1}^{\infty} (-1)^n \mathbf{f}_n(a) t^n \quad (5.3)$$

where

$$\lambda_t(a) = 1 + \sum_{n=1}^{\infty} \lambda_n(a) t^n \tag{5.4}$$

are the ring endomorphisms

$$\mathbf{f}_n : x_j \mapsto x_j^n \tag{5.5}$$

There are well-known determinantal relations between the $\lambda_n$ and the $\mathbf{f}_n$ as follows

$$n!\lambda_n(a) = \det \begin{pmatrix} \mathbf{f}_1(a) & 1 & 0 & \cdots & 0 \\ \mathbf{f}_2(a) & \mathbf{f}_1(a) & 2 & \ddots & \vdots \\ \vdots & \vdots & \ddots & \ddots & 0 \\ \mathbf{f}_{n-1}(a) & \mathbf{f}_{n-2}(a) & \cdots & \mathbf{f}_1(a) & n-1 \\ \mathbf{f}_n(a) & \mathbf{f}_{n-1}(a) & \cdots & \mathbf{f}_2(a) & \mathbf{f}_1(a) \end{pmatrix} \tag{5.6}$$

$$\mathbf{f}_n(a) = \det \begin{pmatrix} \lambda_1(a) & 1 & 0 & \cdots & 0 \\ \lambda_2(a) & \lambda_1(a) & 1 & \ddots & \vdots \\ \vdots & \vdots & \ddots & \ddots & 0 \\ (n-1)\lambda_{n-1}(a) & \lambda_{n-2}(a) & \cdots & \lambda_1(a) & 1 \\ n\lambda_n(a) & \lambda_{n-1}(a) & \cdots & \lambda_2(a) & \lambda_1(a) \end{pmatrix} \tag{5.7}$$

(which come from the Newton relations between the elementary symmetric functions and the power sum symmetric functions).

It follows that the subrings *Symm* and *QSymm* are stable under the $\lambda_n$ and $\mathbf{f}_n$ because $\lambda_n(QSymm \otimes_{\mathbf{Z}} \mathbf{Q}) \subset QSymm \otimes_{\mathbf{Z}} \mathbf{Q}$ by (1.6) and $(QSymm \otimes_{\mathbf{Z}} \mathbf{Q}) \cap \mathbf{Z}[x_1, x_2, \cdots] = QSymm$, and similarly for *Symm*. It follows that *QSymm* and *Symm* have induced $\lambda$ – ring structures. On *Symm* this is of course the standard one.

It follows immediately from (5.5) that

$$\mathbf{f}_n([a_1, a_2, \cdots, a_m]) = [na_1, na_2, \cdots, na_m] \tag{5.8}$$

when $\alpha = [a_1, a_2, \cdots, a_m]$ is seen as a monomial quasisymmetric function as in the realization of *QSymm* described in section 2.

It is more customary to denote the Adams operations (power operations) associated to a $\lambda$-ring structure by $\Psi$'s. However, in the present case they coincide with the Frobenius endomorphisms on *Symm* ([17], section E.2, p.144ff), and their natural extensions to *QSymm*; so it seems natural to use $\mathbf{f}$'s instead in this case.

For any $\lambda$ – ring $R$ there is an associated mapping

$$Symm \times R \longrightarrow R, \quad (\varphi, a) \mapsto \varphi(\lambda_1(a), \lambda_2(a), \cdots, \lambda_n(a), \cdots) \tag{5.9}$$

I.e. write $\varphi \in Symm$ as a polynomial in the elementary symmetric functions $e_1, e_2, \cdots$ and then substitute $\lambda_i(a)$ for $e_i$, $i = 1,2,\cdots$. For a fixed $a \in R$ this is obviously a homomorphism of rings $Symm \longrightarrow R$. We shall often simply write $\varphi(a)$ for $\varphi(\lambda_1(a), \lambda_2(a), \cdots, \lambda_n(a), \cdots)$. Another way to see (5.9) is to observe that for fixed $a \in R$ $(\varphi, a) \mapsto \varphi(\lambda_1(a), \lambda_2(a), \cdots) = \varphi(a)$ is the unique homomorphism of $\lambda$ – rings that takes $e_1$ into $a$. ($Symm$ is the free $\lambda$ – ring on one generator, see also [27].) Note that

$$e_n(\alpha) = \lambda_n(\alpha), \quad p_n(\alpha) = \mathbf{f}_n(\alpha) = [na_1, na_2, \cdots, na_m] \tag{5.10}$$

The first formula of (5.10) is by definition and the second follows from (5.7) because the relations between the $e_n$ and $p_n$ are precisely the same as between the $\lambda_n(a)$ and the $\mathbf{f}_n(a)$.

Now let $P \in NSymm$ be a primitive. Then, by duality

$$\langle P, \alpha\beta \rangle = \langle \mu(P), \alpha \otimes \beta \rangle = \langle 1 \otimes P + P \otimes 1, \alpha \otimes \beta \rangle = 0 \tag{5.11}$$

for any words of length $\geq 1$. Using the Newton relations

$$p_n(\alpha) = p_{n-1}(\alpha)e_1(\alpha) - p_{n-2}(\alpha)e_2(\alpha) + \cdots (-1)^{n-2} p_1(\alpha)e_{n-1}(\alpha) + (-1)^{n-1} ne_n(\alpha) \tag{5.12}$$

it follows that for any primitive in $NSymm$

$$\langle P, p_n(\alpha) \rangle = \pm n \langle P, e_n(\alpha) \rangle \tag{5.13}$$

Now for any $\alpha \in LYN$, $\alpha = [a_1, a_2, \cdots a_m]$ let $\alpha_{red} = [g(\alpha)^{-1}a_1, g(\alpha)^{-1}a_2, \cdots, g(\alpha)^{-1}a_m]$ and define

$$E_\alpha = e_{g(\alpha)}(\alpha_{red}) \tag{5.14}$$

    5.15. *Theorem.* The $E_\alpha$, $\alpha \in LYN$ form a free polynomial basis of $QSymm$ over the integers

For a rather simple direct proof of this (based on Chen - Fox - Lyndon factorization, [4]), see [23].
    Sometimes it is useful to relabel the $E_\alpha$, $\alpha \in LYN$ a bit. Let $eLYN$ be the set of elementary Lyndon words, i.e those Lyndon words $\alpha$ for which $g(\alpha) = 1$. Then the

$$e_n(\alpha), \quad \alpha \in eLYN, \quad n = 1,2,3,\cdots \tag{5.16}$$

are a sometimes convenient relabeling of the free polynomial basis $E_\alpha$, $\alpha \in LYN$. Note that for a fixed $\alpha \in eLYN$, the $e_n(\alpha)$ generate a subalgebra of $QSymm$ that is isomorphic to $Symm$.

Now, for all weights $n$, consider the matrices of integers

$$\left(\langle P_\alpha, E_\beta \rangle\right)_{\alpha, \beta \in LYN_n} \tag{5.17}$$

where $LYN_n$ is the set of Lyndon words of weight $n$ and the columns and rows of (5.17) are ordered by increasing wll-order.
    Claim: this matrix is diagonal with entries $\pm 1$ on the main diagonal. The triangularity follows from theorem 4.24 (ii). As to the diagonal part:

$$g(\alpha)\langle P_\alpha, E_\alpha\rangle = g(\alpha)\langle P_\alpha, \lambda_{g(\alpha)}(\alpha_{red})\rangle \text{ by definition}$$
$$= \pm\langle P_\alpha, p_{g(\alpha)}(\alpha_{red})\rangle \text{ by (5.13)}$$
$$= \pm\langle P_\alpha, \alpha\rangle \text{ by (5.10)}$$
$$= \pm g(\alpha) \text{ by (4.25)}$$

Using that there are precisely $\beta_n = \# LYN_n$ $P_\alpha$'s and $E_\alpha$'s with $\alpha$ of weight $n$, the invertibility (over the integers) of the matrix (5.17) immediately implies both that the $P_\alpha$ are a basis of Prim(*NSymm*) and that the $E_\alpha$ are a free polynomial basis for *QSymm* over the integers

## 6. Frobenius and Verschiebung on *NSymm* and *QSymm*.

To fix notations let

$$\mathbf{f}_n^{Symm} \text{ and } \mathbf{v}_n^{Symm}$$

be the classical Frobenius and Verschiebung Hopf algebra endomorphisms over the integers of *Symm*, characterized by

$$\mathbf{f}_n^{Symm}(p_k) = p_{nk}, \quad \mathbf{v}_n^{Symm}(p_k) = \begin{cases} np_{k/n} & \text{if } n \text{ divides } k \\ 0 & \text{otherwise} \end{cases} \quad (6.1)$$

Now *NSymm* comes with a canonical projection *NSymm* $\longrightarrow$ *Symm*, $Z_n \mapsto h_n$ and there is the inclusion *Symm* $\subset$ *QSymm* (see section 2). The question arises whether ther are lifts $\mathbf{f}_n^{NSymm}$, $\mathbf{v}_n^{NSymm}$ on *NSymm* and extensions $\mathbf{f}_n^{QSymm}$, $\mathbf{v}_n^{QSymm}$ on *QSymm* that satisfy respectively

(i) $\mathbf{f}_n^{?Symm}\mathbf{f}_m^{?Symm} = \mathbf{f}_{nm}^{?Symm}$

(ii) $\mathbf{f}_n$ is homogeneous of degree $n$, i.e. $\mathbf{f}_n(?Symm_k) \subset ?Symm_{nk}$

(iii) $\mathbf{f}_1 = \mathbf{v}_1 = id$

(iv) $\mathbf{f}_m\mathbf{v}_n = \mathbf{v}_n\mathbf{f}_m$ if $(n,m) = 1$ \hfill (6.2)

(v) $\mathbf{v}_n^{?Symm}\mathbf{v}_m^{?Symm} = \mathbf{v}_{nm}^{?Symm}$

(vi) $\mathbf{v}_n$ is homogeneous of degree $n^{-1}$,

i.e. $\mathbf{v}_n(?Symm_k) \subset \begin{cases} ?Symm_{n^{-1}k} & \text{if } n \text{ divides } k \\ 0 & \text{otherwise} \end{cases}$

Here '?' can be 'N', or 'Q'.

Now there exist a natural lifts of the $\mathbf{v}_n^{Symm}$ to *NSymm* given by the Hopf algebra endomorphisms

$$\mathbf{v}_n^{NSymm}(Z_k) = \begin{cases} Z_{k/n} & \text{if } n \text{ divides } k \\ 0 & \text{otherwise} \end{cases} \quad (6.3)$$

and there exist natural extension of the Frobenius morphisms on *Symm* to *QSymm* $\supset$ *Symm* given by the Hopf algebra endomorphisms

$$\mathbf{f}_n^{QSymm}([a_1,\cdots,a_m]) = [na_1,\cdots,na_m] \tag{6.4}$$

which, moreover, have the Frobenius-like property

$$\mathbf{f}_p^{QSymm}(\alpha) = \alpha^p \mod p \tag{6.5}$$

for each prime number $p$. These two families are so natural and beautiful that nothing better can be expected and in the following these are fixed as the Verschiebung morphisms on $NSymm$ and Frobenius morphisms on $QSymm$. They are also dual to each other.

The question to be examined in this section is whether there are supplementary families of morphisms $\mathbf{f}_n$ on $NSymm$, respectively $\mathbf{v}_n$ on $QSymm$, such that (6.2) holds. The first result is negative

6.6. *Theorem.* There are no (Verschiebung-like) coalgebra endomorphisms $\mathbf{v}_n$ of $QSymm$ that extend the $\mathbf{v}_n$ on $Symm$, such that parts (iii)-(vi) of (6.2) hold. Dually there are no (Frobenius-like) algebra homomorphisms of $NSymm$ that lift the $\mathbf{f}_n$ on $Symm$ such that parts (i)-(iv) of (6.2) hold.

It is the coalgebra morphism property which makes it difficult for (6.2) (iii)-(vi) to hold. It is not particularly difficult to find algebra endomorphisms of $QSymm$ that do the job. For instance define the $\mathbf{v}_n$ on $QSymm$ as the algebra endomorphisms given on the generators (5.16) by

$$\mathbf{v}_n(e_k(\alpha)) = \begin{cases} e_{k/n}(\alpha) & \text{if } n \text{ divides } k \\ 0 & \text{otherwise} \end{cases} \tag{6.6}$$

It then follows from (5.12) that

$$\mathbf{v}_n(p_k(\alpha)) = \begin{cases} np_{k/n}(\alpha) & \text{if } n \text{ divides } k \\ 0 & \text{otherwise} \end{cases} \tag{6.7}$$

and all the properties (6.2) follow.

The last topic I would like to discuss is whether there are Hopf algebra endomorphisms of $QSymm$ (and dually, $NSymm$) such that some weaker versions of (6.2) hold.

To this end we first discuss a filtration by Hopf subalgebras of $QSymm$. Define

$$G_i(QSymm) = \sum_{\alpha,\, \lg(\alpha) \leq i} \mathbf{Z}\alpha \subset QSymm \tag{6.8}$$

the free subgroup spanned by all $\alpha$ of length $\leq i$, and let

$$F_i(QSymm) = \mathbf{Z}[e_n(\alpha): \lg(\alpha) \leq i] \tag{6.9}$$

the subalgebra spanned by those generators $e_n(\alpha)$, $\alpha \in eLYN$, $\lg(\alpha) \leq i$ of length less or equal to $i$. Note that this does not mean that the elements of $F_i(QSymm)$ are bounded in length. For instance $F_1(QSymm) = Symm \subset QSymm$ contains the elements

$$\underbrace{[1,1,\cdots,1]}_{n} = e_n = e_n([1]) \tag{6.10}$$

for any $n$.

6.11. *Theorem.* $G_i(QSymm) \subset F_i(QSymm)$

This is a consequence of the proof of the free generation theorem 5.15. Moreover,

$$F_i(QSymm) \otimes_{\mathbf{Z}} \mathbf{Q} = \mathbf{Z}[p_n(\alpha): \lg(\alpha) \leq i]$$
$$F_i(QSymm) = \mathbf{Z}[p_n(\alpha): \lg(\alpha) \leq i] \cap QSymm \qquad (6.11)$$

It follows from (6.11) and (6.12) that the $F_i(QSymm)$ are not only subalgebras but sub Hopf algebras (because $\mu(p_n(\alpha)) = \sum_{\alpha' * \alpha'' = \alpha} p_n(\alpha') \otimes p_n(\alpha''))$.[4]

Now consider a coalgebra endomorphism of $QSymm$. Because of the commutative cofreeness of $QSymm$ as a coalgebra over the module $t\mathbf{Z}[t]$ for the projection

$$QSymm \longrightarrow t\mathbf{Z}[t], [\,] \mapsto 0, [n] \mapsto t^n, \alpha \mapsto 0 \text{ if } \lg(\alpha) \geq 2$$

or, equivalently, because of the freeness of $NSymm$ over its submodule $\sum_{i=1}^{\infty} \mathbf{Z} Z_i$, a homogeneous coalgebra morphism of degree $n^{-1}$ of $QSymm$ is necessarily given by an expression of the form

$$\mathbf{v}_\varphi(\alpha) = \sum_{\alpha_1 * \cdots * \alpha_r = \alpha} \varphi(\alpha_1) \cdots \varphi(\alpha_r)[n^{-1}\text{wt}(\alpha_1), \cdots, n^{-1}\text{wt}(\alpha_r)] \qquad (6.12)$$

for some morphism of Abelian groups $\varphi: QSymm \longrightarrow \mathbf{Z}$. The endomorphism $\mathbf{v}_\varphi$ is a Hopf algebra endomorphism iff $\varphi$ is a morphism of algebras.

6.13. *Proposition.* $\mathbf{v}_\varphi(F_i(QSymm)) \subset F_i(QSymm)$

One particularly interesting family of $\varphi$ 's is the family of ring morphisms given by

$$\tau_n(e_n([1]) = \tau_n(e_n) = (-1)^{n-1}$$
$$\tau_n(e_k(\alpha)) = 0 \text{ for } k \neq n \text{ or } \lg(\alpha) \geq 2 \ (\alpha \in eLYN) \qquad (6.14)$$

Let $\mathbf{v}_n$ be the Verschiebung type Hopf algebra endomorphism defined by $\tau_n$ according to formula (6.12). Then

6.15. *Theorem.*

(i) $\mathbf{v}_n([a_1, \cdots, a_m]) \equiv \begin{cases} n^m[n^{-1}a_1, \cdots n^{-1}a_m] \mod(\text{length } m-1) & \text{if } n \mid a_i \ \forall i \\ 0 \mod(\text{length } m-1) & \text{otherwise} \end{cases}$

(ii) $\mathbf{v}_n$ extends $\mathbf{v}_n^{Symm}$ on $Symm = F_1(QSymm) \subset QSymm$

(iii) $\mathbf{v}_p \mathbf{v}_q = \mathbf{v}_q \mathbf{v}_p$ on $F_2(QSymm)$

(iv) $\mathbf{v}_n \mathbf{f}_n(\alpha) = n^{\lg(\alpha)} \alpha \mod(F_{\lg(\alpha)-1}(QSymm))$

---

[4] I believe the corresponding Hopf ideals in $NSymm$ to be the iterated commutator ideals.

And of course there is a corresponding dual theorem concerning Frobenius type endomorphisms of *NSymm*.

This seems about the best one can do. One unsatisfactory aspect of theorem 6.15 is that there are also other families $\mathbf{v}_n$ that work.

To conclude I would like to conjecture a stronger version of theorem 6.6, viz that there is no family $\mathbf{f}_n$ of algebra endomorphisms of *NSymm* over the integers that satisfies (6.2) (i)-(iii) and such that these $\mathbf{f}_n$ descend to the $\mathbf{f}_n^{Symm}$ on *Symm*.